\documentclass{article}
\usepackage{amsmath,amsfonts,amssymb,amscd}
\usepackage{enumerate}

\usepackage[all]{xy}

\renewcommand{\to}{\ensuremath{\longrightarrow}}

\newcommand{\Z}{\ensuremath{\mathbb Z}}

\renewcommand\theenumi{\@alph\c@enumi}
\renewcommand\theenumii{\@alph\c@enumii}
\renewcommand\theenumiii{\@alph\c@enumiii}
\renewcommand\theenumiv{\@alph\c@enumiv}

\makeatother

\newtheorem{thm}{Theorem}
\newtheorem{lem}[thm]{Lemma}
\newtheorem{prop}[thm]{Proposition}
\newtheorem{cor}[thm]{Corollary}
\newtheorem{defi}[thm]{Definition}
\newtheorem{nota}[thm]{Notation}

\newtheorem{rem}[thm]{Remark}

\newcommand{\eop}{%
  \relax
  \ifvmode
    \noindent
  \else
    \unskip 
    \hskip0pt plus-1fill\relax 
  \fi
  \vrule width0pt 
  \nobreak
  \hfill 
  {\hspace*{\fill}\vrule width3pt height8pt depth0pt}%
}
\newenvironment{proof}{\noindent\emph{Proof.}}
{\eop\par\vspace{\topsep}\vspace{\partopsep}}

\newcommand{\ai}{\vbox to 7pt{\hbox to 7pt{\vrule height 7pt width 7pt}}}

\newcommand{\bs}{\bigskip}

\newcommand{\Spin}{\rm Spin\,}

\newcommand{\qed}{\quad\hfill$\square$}

\begin{document}

\title{Classifications of special double-coverings associated to a non-orientable surface. }

\author{
Anne Bauval\\
Institut de  Math\'ematiques de Toulouse\\
Laboratoire Emile Picard, UMR 5580\\
Universit\'e Toulouse III\\
118 Route de Narbonne, 31062 Toulouse C\'edex - France\\
{\it e-mail: "anne.bauval@math.univ-toulouse.fr"}\\
\\
Claude Hayat\\
Institut de  Math\'ematiques de Toulouse\\
Laboratoire Emile Picard, UMR 5580\\
Universit\'e Toulouse III\\
118 Route de Narbonne, 31062 Toulouse C\'edex - France\\
{\it e-mail: "claude.hayat@math.univ-toulouse.fr"}\\
}

\maketitle

\footnote{57M10; 57R15; 57S05;  57M60. Covering spaces;
 Specialized structures on manifolds (spin manifolds);
  Groups of homeomorphisms or diffeomorphisms;
Group actions in low dimensions
}

\date

\begin{center}
   \begin{minipage}{3.5in}
   \begin{center}
   {\large\it Abstract}
   \end{center}
 This paper  investigates some  possible actions "\`a la Johnson" on the set with $2^{g+1}$ elements, denoted by ${\cal E}$, of  Spin-structures which are interpreted  as  special $\Z_2$-coverings of a trivial $S^1-$fibration  over  a non-orientable surface  $N_{g+1}$. The  group acting is first  a group of orthogonal isomorphisms of  $(H_1(N_{g+1};\Z_2),\cdot)$. A second approach is to consider the subspace of ${\cal E}$ (with $2^{g}$ elements) coming from
special $\Z_2$-coverings of $S^1\times F_g$, where  $F_g\to N_{g+1}$ is the orientation covering. The group acting now is a subgroup of the symplectic group of isomorphisms of the symplectic space $(H_1(F_{g+};\Z_2),\cdot)$. In  both situations, we obtain results on the number of orbits and the number of elements in each orbit. Except in one case, these results do not depend on any necessary choices. We also compare both previous classifications to a third one~: weak-equivalence of coverings
\end{minipage}
 \end{center}

When a surface $N_{g+1}$ is not orientable, the only $S^1-$ principal fibration over $N_{g+1}$ admitting a Spin-structure is the trivial one. We denote it by $P_n$. Like in [BGHM], we interpret a Spin-structure as a special $\Z_2$-covering of the fibration $P_n$. ``Special'' means that the {\sl fibre\/} of $P_n$ is $\Z_2$-covered. The set of these Spin-structures, denoted by ${\cal E}$, is an affine space. It has $2^{g+1}$ elements.

The purpose of this paper is to investigate the action on ${\cal E}$ of a group of homeomorphisms of the surface $N_{g+1}$. We know, for instance from [L], [G], [ZVC], that homeomorphisms of $N_{g+1}$ map onto automorphisms of $(H_1(N_{g+1};\Z_2),\cdot)$, the $\Z_2$-vector space $H_1(N_{g+1};\Z_2)$ endowed with the (non-degenerate) intersection form $\cdot$. The affine structure of ${\cal E}$ compels to fix basis in $\Z_2$-vector spaces. We obtain results on the number of orbits and the number of elements in each orbit. Except in one case, see theorem (\ref{thm:2=4}), these results neither depend on this choice of basis nor on the additional choices in the definition of the action.

The group acting on ${\cal E}$ may be identified with the orthogonal group $O(\Z_2, g)$. Let us remark that this identification depends on the choice of an orthogonal basis of $(H_1(N_{g+1};\Z_2),\cdot)$. The action is an action "\`a la Johnson", i.e. obtained after fixing a section $H_1(N_{g+1};\Z_2)\to H_1(P_n;\Z_2)$, so as to choose a particular way to lift any orthogonal automorphism of $(H_1(N_{g+1};\Z_2),\cdot)$ to an automorphism of $H_1(P_n;\Z_2)$.

A different approach to the classification of special $\Z_2$-coverings of $P_n=S^1\times N_{g+1}$ is to consider the orientation covering $F_g\to N_{g+1}$ and the classification, in [BGHM], of special $\Z_2$-coverings of $P_o=S^1\times F_g$ ($F_g$ is a closed orientable surface). Instead of classifying the $2^{g+1}$ elements of ${\cal E}$, we are thus led to classify only $2^g$ of them, namely those coming from
special $\Z_2$-coverings of $P_o$. This also constrains to restrict the symplectic group which acts.
Surprisingly, although the action we will consider is a restriction (both of the group and of the space on which it acts), it is classified by the same invariant deduced from the classical Arf invariant.

In the third part of this paper, we study the relationship between the group $O(\Z_2,g)$ and the symplectic group. This allows us to compare both previous classifications to a third one~: weak-equivalence of coverings.
A main point is the characterization, among automorphisms of $H_1(P_n;\Z_2)$, of those which are realizable as automorphisms of $P_n$.

We reformulated most classical results on automorphisms of orientable and non-orientable surfaces \cite{GP}, \cite{K2}, \cite{LI}, \cite{ZVC} and on the relationship between $\pi_1F_g$ and $\pi_1N_{g+1}$ \cite{BC} in a convenient way for our purpose.

\section{Special $\Z_2-$coverings}\label{ssec:q}

\bs

We denote by $N_{g+1}$ a non-orientable closed surface and
consider $S^1\hookrightarrow P_n \stackrel{p_n}\longrightarrow N_{g+1}$ a $SO(2)-$principal fibration over $N_{g+1}$ $g>0$. The exact sequence of the fundamental groups is :
$$0\to\Z\to\pi_1P_n\stackrel{(p_n)_\sharp}\longrightarrow\pi_1N_{g+1}\to 0.$$

Let us choose  a presentation of $\pi_1N_{g+1}$ :
\begin{equation}\label{eqi:Ng}
\pi_1N_{g+1}=\langle x_0,\cdots ,x_g\mid \Pi x_i^2\rangle,
\end{equation}
a presentation of $\pi_1 P_n$  given by:
\begin{equation}\label{eqi:Pn}
\pi_1 P_n=\langle u_0,\cdots ,u_g,h\mid  [u_i,h],0\leq i\leq g; \Pi{u_i}^2h^q\rangle.
\end{equation}
where the elements $u_i\in \pi_1P_n$ verify $(p_n)_\sharp(u_i)=x_i$ and  $h$ denotes the image of the generator of $\pi_1S^1$.

 The elements of
$H^1(N_{g+1};\Z_2)$, which is equal to $\Z_2$,  classifie the  $SO(2)-$principal fibrations over $ N_{g+1}. $ Hence there exist only
two such $SO(2)-$principal fibrations depending on the parity of their Chern class : the trivial one and another
one.

\begin{defi}\label{defi:sr}
The  $\Z_2-$covering $E_\psi$  of $P_n$ determined by the exact sequence
$$1\to\pi_1 E_\psi\to\pi_1P_n\stackrel{\psi}\longrightarrow\Z_2\to 0$$
  is called a {\rm special  $\Z_2-$covering} when the homomorphism $\psi$ verifies :
$$ \psi (h)=1\in\Z_2\sim\{0,1\},$$
where $h$ is the image in $\pi_1P_n$ of the generator of $\pi_1S^1=\Z.$
\\
Let us denote by ${\cal E}(P_n)$ the space of special $\Z_2-$coverings $E_\psi.$
\end{defi}

As we know (for examples \cite{A}, \cite{BGHM}), the affine space of $\Spin-$structures   associated to $S^1\hookrightarrow P_n \stackrel{p_n}\longrightarrow N_{g+1}$ coincides with  the space ${\cal E}(P_n)$ of special  $\Z_2-$coverings which may be empty.  The condition $\psi(h)=1$ implies that ${\cal E}(P_n)$
is not empty if and only if $q$ is even which we assume henceforth. That means that the fibration $S^1\hookrightarrow P_n \stackrel{p_n}\longrightarrow N_{g+1}$, we will work on, is trivial.

\bs

\subsection{A presentation of $\pi_1E_\psi=\ker\psi, \psi\in {\cal E}(P_n)$}\label{ssec:kerno}

\bs

 Reidemeister-Schreier's method  [ ZVC] gives a presentation of $\ker\psi$. In the particular case where $S^1\hookrightarrow P_n \stackrel{p_n}\longrightarrow N_{g+1}$ is trivial, i.e.
$\pi_1 P_n=\langle u_0,\cdots ,u_g,h\mid  [u_i,h],0\leq i\leq g; \Pi{u_i}^2\rangle,$ we give a shorter proof.

\begin{thm}\label{thm:kernon}
A presentation of a special  $\Z_2-$covering $E_\psi$  of $P_n$ is given by
$$
\pi_1E_\psi= \langle w_0, \ldots, w_g ,k\mid [w_j,k], \ 0 \leq j \leq g,
                  \prod_{j=0}^g w_j^2 \cdot k^{\Sigma_{i=0}^{i=g}\psi(u_i)} \rangle ,  $$
$u_i$ are the generators of the presentation (\ref{eqi:Pn}) of $\pi_1P_n.$
\end{thm}

\proof \quad For each $i\in \{0,\cdots,g\}$, we choose any $n_i\in\Z$ in the class modulo $2$ of $\psi(u_i).$ We get new generators $w_i=u_ih^{-n_i}$ and a presentation of $\pi_1 P_n$:
$$\pi_1 P_n=\langle w_0,\cdots ,w_g,h\mid  [w_i,h],0\leq i\leq g; \Pi{w_i}^2h^{2\Sigma n_i}\rangle.$$
Let $G=F(x_0,\cdots,x_g)\times F(y)$ be the direct product of the free group generated by $\{x_0,\cdots,x_g\}$ and the free group $F(y)$ generated by one element $y$ and   denote by $H$ the normal subgroup of $G$ generated by
$z=(\Pi x_i^2)y^{2\Sigma n_i}.$ It is clear that $\pi_1 P_n=G/H$.

We denote by $G'$ the direct product of the free group generated by $\{x_0,\cdots,x_g\}$ by the free group $F(y^2)$ generated by the element $y^2$, and  $H'$ the normal subgroup of $G'$ generated by $z=(\Pi x_i^2)y^{2\Sigma n_i}.$  It is clear that $\ker \psi=G'/(H\cap G'),$ and also that $H=H'$ because the set of conjugates of $z$ by elements of $G$ is equal to the set of conjugates of $z$ by elements of $G'$ which  is equal to the set of conjugates of $z$ by elements of $F(x_0,\cdots,x_g).$\\
The conclusion is:
$$\ker \psi=G'/H'.$$

\qed

\bs
\begin{rem}\label{rem:eqf}
\begin{itemize}
\item Contrary to the situation where the surface is orientable (see for example {\rm \cite{BGHM}}), the  fundamental group of the total spaces $E_\psi$ of the  special  $\Z_2-$coverings associated to the trivial bundle over a non-orientable surface $N_{g+1}$ are not all isomorphic. See the theorems (\ref{thm:kernon}) and (\ref{thm:an}).
\item The monomorphism $\pi_1E_\psi\to \pi_1P_n$ is defined by $k \mapsto h^2$  and  $w_i\mapsto u_ih^{\epsilon_i}$ with $\epsilon_i=0$ when $\psi(u_i)=0$ and $\epsilon_i=1$ when $\psi(u_i)=1.$ Hence there are $2^{g+1}$ non-conjugate  special $\Z_2-$ coverings.
\end{itemize}
\end{rem}
\bs

\section{Specific actions of $O(\Z_2^{g+1})$ on ${\cal E}(P_n)$.}\label{sec:1action}

\bs
This action is  an adaptation of Johnson'action \cite{Joh2}. The intersection product defined on $H_1(N_{g+1};\Z_2)$ is now an inner product. We are interested in the group of orthogonal automorphisms of $H_1(N_{g+1};\Z_2)$ because they  are realizable   by homeomorphims of $N_{g+1}$ (see  for example \cite{GP}, \cite{MCP}). We will fix an orthogonal  basis in $H_1(N_{g+1};\Z_2)$ and identify the group of orthogonal automorphisms with the orthogonal group  $O(\Z_2^{g+1})$.
\bs

\subsection{Choice of a lifting of $(p_n)_\star\colon H_1(P_n;\Z_2)\to H_1(N_{g+1};\Z_2)$}\label{ssec:sigma}
\bs

Let us denote by $h_M$ the composition $\pi_1M\to H_1(M;\Z)\to H_1(M;Z_2)$ where the first morphism is the Hurewicz epimorphism and the second one is the change of coefficients. A homomorphism $\psi\colon\pi_1P_n\to\Z_2$  determines a unique $\widetilde\psi\colon H_1(P_n;\Z_2)\to\Z_2$ such that
$\psi=\widetilde\psi\circ h_{P_n}$. This remark allows us to work at the level of homology with coefficients in $\Z_2.$ By abuse of notation, we omit the symbol "tilde" over $\psi$, at least when no confusion is possible.

\bs
\begin{equation}\label{eqi:gene}
\xymatrix{
\pi_1 P_n\ar[dd]_{h_{P_n} }\ar[rr]^{(p_n)_\sharp} &&\pi_1 N_{g+1}\ar[dd]^{h_{N_{g+1}}} \\
&&\\
H_1(P_n;\Z_2)\ar[rr]^{(p_n)_\star} &&H_1(N_{g+1};\Z_2)
}
\end{equation}

We chose $\{x_i\}$ a family of generators of $\pi_1N_{g+1}$ (\ref{eqi:Ng}).
 The family $\{v_i,0\leq i\leq g\}, v_i:=h_{N_{g+1}}(x_i)$ is an orthogonal basis of $H_1(N_{g+1};\Z_2),$ where
the intersection product on $H_1(N_{g+1};\Z_2)$ is

\begin{equation}\label{eqi:oprdct}
H_1(N_{g+1};\Z_2)\times H_1(N_{g+1};\Z_2)\to\Z_2; (v_i,v_j)\mapsto <v_i,v_j>=\delta_{ij}.
\end{equation}

\bs

\begin{nota}\label{nota:sigm}
The fibration $S^1\hookrightarrow P_n\stackrel{p_n}\longrightarrow N_{g+1}$ is trivial hence   $(p_n)_\star$ admits {\sl linear sections}. We choose arbitrarily one of them, it means that we choose an embedding of $\pi_1N_{g+1}\hookrightarrow\pi_1P_n,$ and this choice will induce one for
$H_1(N_{g+1};\Z_2)\hookrightarrow H_1(P_n;\Z_2); x\mapsto \overline x.$ 

Since (by definition) $v_i=h_{N_{g+1}}(x_i)$, we get $\overline v_i=h_{P_n}(u_i)$, hence $\{\overline v_i,h\}$ is a basis of $H_1(P_n;\Z_2).$ A linear section $\sigma$ of $(p_n)_\star$ is determined by its value on the basis $\{v_i\}$ of $H_1(N_{g+1};\Z_2)$ by $\sigma(v_i)=\overline v_i+\rho_ih$, with $\rho_i\in \Z_2,$ arbitrarily fixed.
\end{nota}

\bs

\begin{nota}\label{nota:Fsig}
Let $F\colon H_1(N_{g+1};\Z_2)\to H_1(N_{g+1};\Z_2)$ be an  automorphism, a lifting
$$F_\sigma\colon H_1(P_n;\Z_2)\to H_1(P_n;\Z_2)$$ is defined by the relations:
$$F_\sigma(\sigma (v_i))=\sigma (F(v_i)), 0\leq i\leq g; \quad F_{\sigma} (h)=h,$$
and extended by linearity.
\end{nota}

\bs


If $A=(a_{ij})$ is the matrix of $F$ in the basis $\{v_i\}$ then the matrix of $F_\sigma$ in the basis $\{\overline v_i,h\}$ is
$\left( \begin{array}{cc}A&0\\ d&1  \end{array}\right)$ where $d$ is a line with
 $d_j=\Sigma \rho_ia_{ij}+\rho_j.$

\begin{prop}\label{prop:J_n}
The map $J_n\colon O(\Z_2^{g+1})\to SL(\Z_2^{g+2})$ defined as above by $J_n(A)=\left( \begin{array}{cc}A&0\\ d&1  \end{array}\right) $ is a monomorphism.
\end{prop}

\qed

\bs

 \begin{defi}\label{defi:1} A right  action $A_1^\sigma$ of $O(\Z_2^{g+1})$ on ${\cal E}(P_n)$  is defined by :
$$(E_\psi,F)\mapsto E_{\psi\circ F_{\sigma}}.$$
\end{defi}

\bs

\subsection{Classification of linear forms on $\Z_2^{g+1}$ under the action of $O(\Z_2^{g+1})$}\label{ssec:linform}

\bs
The family $\{\sigma(v_i), 0\leq i\leq g,h\}$ is a basis of $H_1(P_n;\Z_2)$. Hence there is a bijection between the special  $\Z_2-$coverings  $E_\psi$ and the  linear forms $\theta=\psi\circ \sigma$ on $\Z_2^{g+1}$ and
 the action $A_1^\sigma$ is determined by the action by right composition of $ O(\Z_2^{g+1})$ on the linear forms $\theta$ on $\Z_2^{g+1}.$

We consider the right action by composition of $O(\Z_2^{g+1})$ on  $\Z_2^{g+1}.$
If $\theta$ and $\theta'$, two linear forms on $\Z_2^{g+1}$, are  in the same orbit,  we will write $\theta\sim\theta'.$

Let us denote by $v=\Sigma_{i=0}^{i=g} v_i$. Any $T\in
O(\Z_2^{g+1})$ lets invariant $v^\perp=\{ x\in \Z_2^{g+1}\mid
0=<x,x>=\Sigma x_i=<v,x>\}$  then also $v^{\perp\perp}=\Z_2 v$ {\sl
i.e.} $T(v)=v$. This is the proof of the following  lemma which will
be useful in this paper.

\begin{lem}\label{lem:uti} Any $T\in O(\Z_2^{g+1})$ lets fixed $v=\Sigma_{i=0}^{i=g} v_i$.
\end{lem}

\qed

\begin{prop}\label{prop:mi} Under the action by right composition by elements of $O(\Z_2^{g+1})$,
a) the following two linear forms are fixed : $\theta_0\colon x\mapsto 0$ and $\theta_1\colon x=\Sigma_{i=0}^{i=g} x_iv_i\mapsto
\Sigma_{i=0}^{i=g} x_i$;

Apart from these two linear forms,

\noindent b) when $g=1$, the two linear forms $x_0v_0+x_1v_1\mapsto x_0$ and $x_0v_0+x_1v_1\mapsto x_1$ are equivalent;

\noindent c) when $g\geq 2$, $\theta\sim\theta'$ if and only if
$\Sigma_{i=0}^{i=g} \theta(v_i)=\Sigma_{i=0}^{i=g} \theta'(v_i)$.
\end{prop}

\bs

\proof

Let us denote by $m_i=\theta(v_i)$ and $m'_i=\theta'(v_i)$.

1) First we want  to show that if $\theta\sim\theta'$ then $\Sigma m_i=\Sigma m'_i$. If $\theta\sim\theta'$ there exists $T$  such that $\theta'=\theta\circ T$, we get $\Sigma m'_i=\theta'(v)=\theta\circ T(v)=\theta(v)=\Sigma m_i.$ because $T(v)=v.$

2) Choose any $a=\Sigma a_iv_i$ with $\Sigma a_i=0$, and define, for $x\in \Z_2^{g+1}$, $T_a(x)=x+<a,x>a$. Here we have the equality $<a,a>=\Sigma a_i^2=\Sigma a_i=0$, so $<T_a(x),T_a(y)>=<x,y>$ and $T_a\in O(\Z_2^{g+1})$.

Let $\theta$ be a linear form such that there exists at least one
pair $i,j$ such that  $m_i\neq m_j.$  Notice that the order of the
terms $0$ and $1$ is not important because the permutations are
elements of $O(\Z_2^{g+1})$ so  we can work with
$m=(0,1,m_2,\cdots,m_{g})$.  Let $a=(1, 1+\lambda, m_2,\cdots,m_g)$
with $\lambda=\Sigma_{i\geq 2}m_i$, we have $\Sigma a_i=0$ and
$\Sigma a_im_i=1$ hence $m'_i=\theta\circ
T_a(v_i)=\theta(v_i+a_ia)=m_i+a_i$ and $m'=(1,\lambda,0,\cdots,0).$

\qed

\bs By definition,  $\Sigma_{i=0}^{i=g} \theta(v_i)\psi\circ\sigma(v_i)=\Sigma_0^g
\psi(\overline v_i)+\Sigma_0^g\rho_i, $ hence we have the following
lemma (for the notation $\overline v_i$ see (\ref{nota:sigm}) :

\begin{lem}\label{lem:01} Let us distinguish $\psi_0\circ\sigma\colon x\mapsto 0$ and $\psi_1\circ\sigma\colon x=\Sigma x_iv_i\mapsto \Sigma x_i, $ (with $\psi_i(h)=1$).
\begin{itemize}
\item  $\Sigma_0^g \psi_0(\bar v_i)=\Sigma_0^g\rho_i;$
\item  $\Sigma_0^g \psi_1(\bar v_i)=g+1+ \Sigma_0^g\rho_i.$ 
\end{itemize}
\end{lem}

A  direct consequence of the above classification is the following theorem :

\begin{thm}\label{thm:1} Let us distinguish two  special $\Z_2-$coverings $ E_{\psi_0}$ and $ E_{\psi_1}$, characterized by the property : $\psi_0\circ\sigma\colon x\mapsto 0$ and $\psi_1\circ\sigma\colon x=\Sigma x_iv_i\mapsto \Sigma x_i.$

1) If $g=1$, under the action $A_1^\sigma$ on the space ${\cal E}$ of special  $\Z_2-$coverings, there are two fixed points
$ E_{\psi_0}$ and $ E_{\psi_1}$. The two other elements are in the same class.

2) If $g\geq 2$, there are also two fixed points
$ E_{\psi_0}$ and $ E_{\psi_1}$ and two other classes
\begin{itemize}
\item if $g$ is even, $ \{E_\psi\mid \Sigma_0^g (\psi(\overline v_i)+\rho_i)=1,\psi\neq\psi_1\},$
$\{E_\psi\mid \Sigma_0^g (\psi(\overline v_i)+\rho_i)=0,\psi\neq \psi_0\};$
\item if $g$ is odd, $ \{E_\psi\mid \Sigma_0^g (\psi(\overline v_i)+\rho_i)=1\},$
$\{E_\psi\mid \Sigma_0^g (\psi(\overline v_i)+\rho_i)=0,\psi\neq \psi_0,\psi_1\};$
\end{itemize}
\end{thm}

\qed

\bs

\subsection{Isotropy subgroups of $O(\Z_2^{g+1})$}\label{ssec:isosg}

\bs

From the proposition (\ref{prop:mi}), we know that, when $g>1$, the
set of linear form on $\Z_2^{g+1}$ is the union of four orbits under
the action of $O(\Z_2^{g+1})$, two of them are reduced to an element
$\{\theta_0\},\{\theta_1\}$ and the two others are $
orb_0=\{\theta\neq\theta_0,\theta_1\mid \Sigma_0^g\theta(v_i)=0\}$
and $ orb_1=\{\theta\neq\theta_0,\theta_1\mid
\Sigma_0^g\theta(v_i)=1\}.$ Moreover any $F\in O(\Z_2^{g+1})$ allways fixes
 $\theta_0$ and $\theta_1$ but possibly no other linear
form (for example: $F$ with matrix $\left(\begin{array}{cc} 0&1\\
1&0\end{array}\right)$ in the basis $(v_0,v_1)$).

\bs
\begin{thm}\label{thm:1234}
\begin{description}
\item[1)] The isotropy subgroups $Iso_{\theta_0}$ of $\theta_0$, and $Iso_{\theta_1}$ of $\theta_1$, are the whole
$O(\Z_2^{g+1})$.
\item[2)] Let us chose  a representative $\alpha_i\in {\rm  orb}_i, (i=0,1)$, for example such that $\alpha_0(v_0)=1;\alpha_0(v_1)=1;\alpha_0(v_j)=0$ for $j\neq 0;1$ and $\alpha_1(v_0)=1,\alpha_1(v_j)=0$ for $j\neq 0.$  For every $\theta\in {\rm orb}_i$ the isotropy subgroup $Iso_\theta$ is conjugate to the isotropic subgroup $Iso_{\alpha_i}$ and
$$Iso_{\alpha_1}=\{F\in O(\Z_2^{g+1})\mid F(v_0)=v_0\}\simeq O(\Z_2^{g});$$
$$Iso_{\alpha_0}=\{F\in O(\Z_2^{g+1})\mid F(v_0+v_1)=v_0+v_1\}.$$
\item[3)] The elements of $Iso_{\alpha_1}$ are products of \\
$\bullet$  permutations of the elements $v_i, 1\leq i\leq g$ and \\
$\bullet$ if $g\geq 4,$ the transvection $T_{v_1+v_2+v_3+v_4}\colon x\mapsto x+<x,v_1+v_2+v_3+v_4>(v_1+v_2+v_3+v_4).$
\item[4)] Here $g\geq 1$. Every element of $Iso_{\alpha_0}$ is a product of\\
$\bullet$ elements of $O(\Z_2^2)\times O(\Z_2^{g-1})$, and\\
$\bullet$ if $g\geq 3$, the transvection $T_{v_0+v_1+v_2+v_3}\colon x\mapsto x+<x,v_0+v_1+v_2+v_3>(v_0+v_1+v_2+v_3).$
\end{description}
\end{thm}
\bs

To prove this theorem, we first establish the following lemma :

\begin{lem}\label{lem:trans}
Let us denote by $H=\{x\in V\mid x\neq 0, x\neq e\}$ where $V=\oplus_{i=1}^{i=g}\Z_2v_i, e=v_1+\cdots+v_g.$ The action of $O(\Z_2^{g})$ on the subset $H_0=\{x\in H\mid <x,e>=0\}$ and on $H_1=\{x\in H\mid <x,e>=1\}$, is transitive.
\end{lem}

 {\sl Proof of the lemma.}

 First we remark that $O(\Z_2^{g})$ contains all the permutations of the elements $v_i, 1\leq i\leq g$ and all the transvections $T_u(x)=x+<x,u>u$ such that $<u,e>=0.$ These transvections verify $T_u^{-1}=T_u.$\\
$i)$ When $x$ is an element of $H_0$, it is possible to suppose (up to permutations) that $<x,v_1>=1$ and $<x,v_2>=0.$ For $u=x+v_1+v_2$, we have $<u,e>=0$ and $T_u(v_1+v_2)=x.$\\
$ii)$ When $x$ is an element of $H_1$,  it is possible to suppose (up to permutations) that $<x,v_1>=0.$  For $u=x+v_1$, we have $<u,e>=0$ and $T_u(v_1)=x.$

\qed

{\sl Proof of the part $3)$ of the theorem :}

Let $V=\oplus_{1\leq i\leq g}\Z_2v_i$, then $Iso_{\alpha_1}$ is naturally isomorphic to $O(V)$ (by restriction). If $g=0$ or $g=1$ then $O(V)=\{{\rm id}_V\}$ and the result is obvious.

Assume henceforth that $g\geq 2$ and that the result is true for $g-1$. The elements of $Iso_{\alpha_1}$ fixing $v_1$ are called of ``type $R$''. By induction hypothesis they are products of the described generators.

Let $\gamma\in O(V)$ and $e=v_1+\ldots+v_g$. By lemma (\ref{lem:uti}), $\gamma(e)=e$. Since $g\geq 2$, we have $v_1\neq e$, hence $\gamma(v_1)\neq\gamma(e)=e$, which means that $\gamma(v_1)$ has at least one coordinate equal to $0$. Therefore (up to permutations) we may assume $\gamma(v_1)\in V':=\oplus_{2\leq i\leq g}\Z_2v_i$.

First case: if $\gamma(v_1)\neq v_2+\ldots+v_g$ or if $\gamma(v_1)=v_2$, the lemma (\ref{lem:trans}) applied to $V'$ says that there exists $F$ of type $R$ such that $F^{-1}(\gamma(v_1))=v_2$, hence $F^{-1}\circ\gamma$ is the product of a permutation and of a type $R$.

Second case: if $\gamma(v_1)=v_2+\ldots+v_g$ and $g>2$ (which implies that $g$ is even and $>2$, hence $g\geq 4$), then $T^{-1}_{v_1+v_2+v_3+v_4}(\gamma(v_1))=v_1+v_5+\ldots+v_g\neq v_2+\ldots+v_g$, hence $T^{-1}_{v_1+v_2+v_3+v_4}\circ\gamma$ belongs to the first case above.

\bs

 {\sl Proof of the part $4)$ of the theorem:}

For $g=1$ it is evident.\\
Let us suppose that $g>1$. For $\gamma\in Iso_{\alpha_0} $ i.e. $\gamma(v_0+v_1)=v_0+v_1$, we write $\gamma(v_0)=x+y$ where $x\in V(v_0;v_1)$, the space generated by $\{v_0,v_1\},$ and $y\in V(v_2;\cdots;v_g)$ the space generated by $\{v_2,\cdots,v_g\}.$
The following computations show that $<y,y>=0$ and that $x=v_0$ or $x=v_1$ hence $<x,x>=1$ :
$$\begin{array}{rl} 1=&<\gamma(v_0),\gamma(v_0>=<x+y,x+y>=<x,x>+<y,y>; \\
0=&<\gamma(v_0),\gamma(v_1)>=<\gamma(v_0),\gamma(v_0)+v_0+v_1>\\
=&<\gamma(v_0),\gamma(v_0)>+
<x+y,v_0+v_1>=1+<x,v_0+v_1>.\end{array}
$$

First case: if $y=0.$ then $\gamma\in O(\Z_2^2)\times O(\Z_2^{g-1}).$

\bs

Second case: if $y\neq 0$ (which implies that $g\geq 3$) but $y\neq v_2+\cdots+v_g$, the lemma (\ref{lem:trans}) shows that there exists $F\in id_{\Z_2}\times O(\Z_2^{g-1})$ such that $F^{-1}(y)=v_2+v_3$ then
$$T_{v_0+\cdots+v_3}^{-1}(F^{-1}(\gamma(v_0)))=T_{v_0+\cdots+v_3}(x+v_2+v_3)=v_0+v_1+x,$$
so
$$T_{v_0+\cdots+v_3}^{-1}(F^{-1}(\gamma(v_0)))=\left\{\begin{array}{c} v_0\mbox{ if } x=v_1\\ v_1 \mbox{ if } x=v_0,\end{array}\right\},$$
{\sl i.e.}  we are now in the situation of the first case above.

\bs

Third case : if $y\neq 0$ (which implies that $g\geq 3$) but $y= v_2+\cdots+v_g$, then
$$T_{v_0+\cdots+v_3}^{-1}(\gamma(v_0))=T_{v_0+\cdots+v_3}(x+v_2+\cdots+v_g)=(v_0+v_1+x)+(v_4+\cdots+v_g)=x'+y'.$$
If $g=3$, then $y'=0$ and we are in the first case. If $g\geq 4$, ($y'\neq 0$ but $y'\neq v_2+\cdots+v_g$) we are in the second case.

\qed

\bs

\section{Quadratic action  }\label{sec:a2}
\bs
In this section, we will take advantage of the existence of an oriented double covering $F_g\stackrel{\pi}\longrightarrow N_{g+1}$ of $ N_{g+1}$ and the fact that $\pi_1F_g$ is characteristic in $\pi_1 N_{g+1}.$ Then it is possible  to use the classification obtained for special  $\Z_2-$coverings associated to a
$SO(2)-$principal bundle over $F_g$ [BGHM].

\bs
\subsection{Oriented double covering}\label{ssec:odc}

\bs

Let $S^1\hookrightarrow P_o\stackrel{p_o}\longrightarrow F_{g}$ be  the pull-back by $\pi\colon F_g\to N_{g+1}$ of the  fibration $S^1\hookrightarrow P_n \stackrel{p_n}\longrightarrow N_{g+1}$ studied in the first part of the article. We get
the following commutative diagram:
\begin{equation}\label{eqi:diagm1}
\xymatrix{
S^1\ar[dd] \ar[rr]&&P_o\ar[dd]_{\widetilde\pi }\ar[rr]^{p_o} &&F_g\ar[dd]_{\pi} \\
&&&&\\
S^1\ar[rr]&&P_n\ar[rr]^{p_n} &&N_{g+1}
}
\end{equation}
Moreover, since $S^1\hookrightarrow P_n \stackrel{p_n}\longrightarrow N_{g+1}$ was trivial $S^1\hookrightarrow P_o\stackrel{p_o}\longrightarrow F_{g}$ is also trivial. More precisely, any trivialisation of the former induces a trivialisation of the latter, compatible with the maps $\pi$ and $\tilde\pi.$ Thus, if we still denote by $y\mapsto \overline y; H_1(F_g;\Z_2)\to H_1(P_o;\Z_2)$, the linear section of $(p_o)_\star$ induced by the section $x\mapsto \overline x; H_1(N_{g+1};\Z_2)\to H_1(P_n;\Z_2)$ of $(p_n)_\star$ chosen in (\ref{nota:sigm}), we get $\tilde\pi_\star(\overline y)=\overline{\pi_\star}(y).$

\bs 

Reidmeister-Schreier's method (for example as in [ZVC]) furnishes a symplectic basis $(c_i)_{1\leq i\leq 2g}$ of $(H_1(F_g;\Z_2),\dot)$ ({\sl i.e.} all the intersection products are $0$ except $(c_{2i-1},c_{2i})=1, 1\leq i\leq g$) such that $\pi_\star(c_{2i})=\pi_\star(c_{2i-1})=v_0+v_i.$  In particular,  the sections defined in the above paragraph verify $\tilde\pi(\overline{c_{2i-1}})=\tilde\pi(\overline{c_{2i}})=\overline{v_0+v_i}.$

\begin{equation}\label{eqi:overline}
\xymatrix{
\pi_1F_g\ar[dd]_{h_{F_g}}\ar[rr]^{\pi_\sharp} &&\pi_1N_{g+1}\ar[dd]^{h_{N_{g+1}}} \\
&&\\
H_1(F_g;\Z_2)\ar[rr]_{\pi_\star}&&H_1(N_{g+1};\Z_2)\\
&&\\
H_1(P_o;\Z_2)\ar[uu]^{(p_o)_\star}\ar[rr]_{\widetilde\pi_\star}\ar[dr]_{\psi}&&H_1(P_n;\Z_2)\ar[uu]_{(p_n)_\star}\ar[dl]^{\varphi}\\
 &\Z_2}
\end{equation}

In this diagram and in the following, the map $h_M$ is (like in subsection \ref{ssec:sigma}) the composition $\pi_1M\to H_1(M;\Z)\to H_1(M;Z_2)$ where the first morphism is the Hurewicz epimorphism.

\begin{nota} \label{nota:epi} ${\cal E}(P_o)$ is the notation for the set of special  $\Z_2-$coverings associated to the principal bundle $P_o$ over $F_g.$ 

The elements of ${\cal E}_\pi$ are the special  $\Z_2-$coverings $E_\varphi$ such that  there exists $\psi\colon H_1(P_n;\Z_2) \to \Z_2$ with $\varphi=\psi\circ \tilde\pi_\star.$
\end{nota}

\bs

\begin{prop} \label{prop:2g} The number of elements of  ${\cal E}_\pi$ is equal to $2^g$.
\end{prop}

\proof\quad  Let us suppose that $\psi_1\circ \tilde\pi_\star=\psi_2\circ \tilde\pi_\star$ then, depending of the values of $\psi(\overline v_0)$ and  $\psi'(\overline v_0)$, there are two possibilities  $$\left\{
\begin{array}{l}
\forall i\in 0\cdots g,  \psi_1(\overline v_i)=\psi_2(\overline v_i)\\
 \mbox{or}  \\
\forall i\in 0\cdots g,  \psi_1(\overline v_i)=\psi_2(\overline v_i)+1
\end{array}
\right.$$
In the second situation, we will write $\psi_1=\psi_2+ {\bf 1}.$\qed 

\bigskip

We recall the classification obtained for special $\Z_2-$coverings associated to a
$SO(2)-$principal bundle over $F_g$ [BGHM].

\bs

\subsection{Choice of a quadratic section of $(p_o)_\star\colon H_1(P_o;\Z_2)\to H_1(F_g;\Z_2)$}\label{ssec:s}

\bs

\begin{prop}\label{prop:s}{\rm [BGHM]}
For 	any family $\{d_i,1\leq i\leq 2g\}$ of $H_1(P_o;\Z_2)$ such that 
$(p_o)_\star(d_i)=c_i$ (which may always be written $d_i=\overline c_i+r_ih$), 
\begin{itemize}
\item $\{d_i,1\leq i\leq 2g\}\cup \{ h\}$ is a  basis of $H_1(P_o;\Z_2);$
\item
there exists an
unique  map $s\colon  H_1(F_g;\Z_2)\to H_1(P_o;\Z_2)$such that  $s(c_i)=d_i,$ and
  \begin{eqnarray}\label{eq:*}
  s(a+b)=s(a)+s(b)+(a.b)h,
\end{eqnarray}
where $a.b$ is the intersection product of any elements $a,b$ of $H_1(F_g;\Z_2);$
\item such a map, called therefore a {\rm quadratic section}, satisfies  $(p_o)_\star\circ s=id.$
\end{itemize}
\end{prop}

\qed

\bs
\begin{defi}\label{defi:widetildef} Let  $L$ be the sympletic
matrix $(a_{ij})_{i,j\leq 2g}$ writen in the basis $\{c_i\}$, of a symplectic
automorphism $f\colon H_1(F_g;\Z_2)\to H_1(F_g;\Z_2)$. We define\\
 $ f_s\colon H_1(P_o;\Z_2)\to H_1(P_o;\Z_2)$ by linearity 
from
$$f_s (s(c_i)):=s(f(c_i)),  f_s(h):=h.$$ 
The quadratic section $s$ is choosen like in Proposition (\rm{\ref{prop:s}}).
\noindent The matrix of $ f_s$ in the basis $\{\overline c_i,h\}$ is:
$\left( \begin{array}{cc}
                                  L&0\\
                                \Delta&1
                                 \end{array}
                                  \right)$
where $\Delta$ is a line with $2g$ terms $b_j=\Sigma a_{ij}r_i+S_j+r_j, S_j=\Sigma a_{2i,j}a_{2i-1,j}. $
\end{defi}

\bs
 A basis is chosen  in $H_1(F_g;\Z_2)$ endowed with the product intersection which is symplectic,  so it is possible to identify the group of symplectic automorphisms of $H_1(F_g;\Z_2)$ with the symplectic group  $Sp(\Z_2,2g)$.

As in [BGHM], the action of $Sp(\Z_2,2g)$ on ${\cal E}(P_o)$ is defined by $(\varphi,f)\mapsto \varphi\circ f_s$.

In order to apply  this action to the subset ${\cal E}_\pi$ (\ref{nota:epi}), we need to restrict to a subgroup of $Sp(\Z_2,2g)$ leaving ${\cal E}_\pi$ stable. 
\bs

Let us consider the subgroup  $G_s=
\{f\in Sp(\Z_2,2g)\ |\ \forall\varphi\in{\cal E}_\pi, \varphi\circ f_s\in{\cal E}_\pi\}.$

\bs

\begin{prop}\label{prop:gs} $G_s=
\{f\in Sp(\Z_2,2g)\ |\ \forall\varphi\in{\cal E}_\pi, \varphi\circ f_s\in{\cal E}_\pi\}=$

$\{f\in Sp(\Z_2,2g)\ |\ f_s(\ker(\tilde\pi_*))\subset \ker(\tilde\pi_*)\}$

\end{prop}

\proof\quad
Using the following equivalences where $u,v,w$ are linear maps :
$$u(\ker v)\subset\ker w\Leftrightarrow wu(\ker v)=0\Leftrightarrow\ker v\subset\ker(wu),$$
and the following lemma,
 we get the proposition.

\qed

\bs
\begin{lem}\label{lem:bot}
$${\cal E}_\pi=\{\varphi\in {\cal E}(P_o)\mid \varphi(\ker\pi)=0\},$$
and 
$$\cap_{\varphi\in {\cal E}_\pi}\ker\varphi=\ker \tilde\pi.$$
\end{lem}

\bs

\bs
\subsection{Classification of the elements of ${\cal E}_{\pi}$}\label{ssec:cale}

\bs
\begin{nota}\label{nota:beta} 
 Let $e_i=c_{2i-1}$, $e'_i=c_{2i-1}+c_{2i}$, and $t=\sum_{1\leq i\leq g}\beta_ie_i$ with $\beta_i=r_{2i-1}+r_{2i}+1$ (see Proposition (\ref{prop:s})). Note that $(e=(e_i)_{1\leq i\leq g},e'=(e'_i)_{1\leq i\leq g})$ is a symplectic basis and $e'$ is a basis of $\ker\pi_*$ (hence $(\ker\pi_*)^\perp=\ker\pi_*$).

 Note that $(\overline e=(\overline e_i)_{1\leq i\leq g},\overline e'=(\overline e'_i)_{1\leq i\leq g}, h)$ is a  basis of $H_1(P_o;\Z_2)$ and $\overline e'$ is a basis of $\ker\tilde\pi_*$ \rm{(see Notations (\ref{nota:sigm}))}.

\end{nota}

In the following diagram, the homology groups are with coefficient in $\Z_2.$

\begin{equation}\label{eqi:diagm2}
\xymatrix{
H_1(P_o)\ar[dd]_{\widetilde\pi_\star }\ar[rr]^{(p_o)_\star} &&H_1(F_g)\ar[dd]_{\pi_\star}\ar[rr]^{f}&&H_1(F_g)\ar[dd]_{\pi_\star}&&H_1(P_o)\ar[dd]_{\widetilde\pi_\star }\ar[ll]_{(p_o)_\star} \\
&&&&&&\\
H_1(P_n)\ar[rr]^{(p_n)_\star} &&H_1(N_{g+1})\ar[rr]^{F}&&H_1(N_{g+1})&&H_1(P_n)\ar[ll]_{(p_n)_\star} 
}
\end{equation}

\bs

In $H_1(F_g,\Z_2)$, we will work with the new symplectic basis $\{ e_i,e'_i\, 1\leq i\leq g\},$ with $e_i=c_{2i-1}; e_i'=c_{2i-1}+c_{2i},$ and in $H_1(N_{g+1};\Z_2)$ with the basis $(v_0,v_0+v_1,\ldots,v_0+v_g).$

The morphism $\pi_\star\colon  H_1(F_g,\Z_2)\to  H_1(N_{g+1},\Z_2),e_i\mapsto v_0+v_i,e'_i\mapsto 0,$ has a kernel $\ker \pi_\star=\oplus_{i=1}^g\Z_2e'_i$ and $Im(\pi_\star)=\oplus_{i=1}^g\Z_2(v_0+v_i).$

The study of the action of $G_s$  (\ref{prop:gs}) on ${\cal E}_\pi$ will naturally lead us to consider the following subgroup $K_t, t=\sum_{1\leq i\leq g}\beta_ie_i$ (see notations \ref{nota:beta}) of $Sp(\Z_2,2g)$, which contains the particular transvections occurring in the proof of the theorem (\ref{thm:2=4}) below. In the proposition (\ref{prop:genKt}) we shall prove that $K_t$ is in fact the smallest subgroup containing these transvections.

\begin{nota}\label{nota:kt}
$K_t=\{f\in Sp(\Z_2,2g)\ |\ {\rm Im}(f-id)\perp\ker\pi_*\oplus\Z_2t\}$.
\end{nota}

\bs

\begin{rem}\label{rem:k0}
$K_t\subset K_0$, and $K_0=\{f\in Sp(\Z_2,2g)\ |\ \pi_*\circ f=\pi_*\}$ (since $\ker\pi_*^\perp=\ker\pi_*$).
\end{rem}

\bs 

The next proposition is a generalization of the well-known case $V=\{0\}$. When applied to $V=\ker\pi_*\oplus\Z_2t$, it shows that the transvections $T_Y\colon x\mapsto T_Y(x)=x+<x,Y>Y$ such that $Y\in\ker\pi_*$ and $Y.t=0$ are generators of $K_t$.

\begin{prop}\label{prop:genKt}
Let $V$ be any linear subspace of $\Z_2^{2g}$. Then the subgroup
$$G_V:=\{ f\in Sp(\Z_2,2g)\ |\ {\rm Im}(f-id)\perp V\}$$
is generated by its transvections.
\end{prop}

\bs

\proof\quad

We shall prove that any $f\in G_V$ is a product of transvections $T_Y\in G_V$, by induction on the dimension of $W_f:={\rm Im}(f-id)$ (for this we shall use twice that for all $Y\in W_f$, we have $T_Y\in G_V$ and $W_{T_Yf}\subset W_f$). When $W_f=\{0\}$, there is nothing to prove. Assume now that $f\neq id$ and that the property holds for all $h\in G_V$ such that ${\rm dim}(W_h)<{\rm dim}(W_f)$.\\

First case. If there exists some $z\in\Z_2^{2g}$ such that $f(z).z=1$, let $Y=f(z)-z$ and $h=T_Y\circ f$. Then $T_Y\in G_V$ and $W_h\subset W_f$, and $W_h\neq W_f$ (since $Y\in W_f$ but for all $x\in\Z_2^g, (h(x)-x).z=0\neq 1=Y.z$). By induction hypothesis, $h$ is a product of transvections belonging to $G_V$, hence so is $f$.

Second case. If for all $z\in\Z_2^{2g},f(z).z=0$, note that for all $u,v\in\Z_2^{2g},
u.f(v)=f(u).v$ (for this, compute $f(u+v).(u+v)$).
Choose any $u\in\Z_2^{2g}$ such that $e:=f(u)-u\neq 0$ and let $k=T_e\circ f$. Then $T_e\in G_V$ and $W_k\subset W_f$. Moreover, for any $v\in\Z_2^{2g}$ such that $e.v=1$, we have $k(v).v=e.f(v)=f(u).f(v)-u.f(v)=u.v-f(u).v=e.v=1$. By the first case, $k$ is a product of transvections belonging to $G_V$, hence so is $f$.\qed

\bs

In [BGHM],  for any $SO(2)-$principal bundle over the orientable $F_g$, once a quadratic section $s$ is chosen like in the proposition (\ref{prop:s}), the autors defined an action of $Sp(\Z_2,2g)$  on the set of  special  $\Z_2-$coverings. They also proved that two special  $\Z_2-$coverings determined by $\varphi$ and $\varphi'$ are in the same orbit under this action il and only if the quadratic forms $\varphi\circ s$ and $\varphi'\circ s$ have the same Arf-invariants denoted here by $\alpha(\varphi\circ s)=\sum_{1\leq i\leq g}\varphi\circ s(e_i)\varphi\circ s(e'_i)$. Surprisingly, although the action we consider here is a restriction (both of the group and of the space on which it acts), it is classified by the same invariant.

\bs

\begin{thm}\label{thm:2=4}
\begin{description}
\item[i)] $K_t\subset G_s$
\item[ii)]  For any $\varphi\in{\cal E}_\pi$, $\alpha(\varphi\circ s)=r+\varphi(\overline t)$, where $r=\sum_{1\leq i\leq g}r_{2i-1}r_{2i}$.
\item[iii)] The actions of $G_s$ and $K_t$ on ${\cal E}_\pi$ have the same orbits, which are classified by the Arf invariant $\alpha(\varphi\circ s).$
\item[iv)] In general there are two orbits, each with $2^{g-1}$ elements. One exception arises with a particular choice of the quadratic section (proposition (\ref{prop:s})): when $r_{2i-1}+r_{2i}=1$ for all indices $i$. In this case, all the elements of 
${\cal E}_\pi$ are in the same orbit.
\end{description}
\end{thm}

\proof

\begin{itemize}
\item iv) is a straightforward consequence of ii) and iii).

\item ii) is mere calculus : 
$$\begin{array}{rl}
\alpha(\varphi\circ s)=&\sum_{1\leq i\leq g}\varphi\circ s(e_i)\varphi\circ s(e'_i)\\
=&
\sum_{1\leq i\leq g}\varphi(\overline e_i+r_{2i-1}h)\varphi(\overline e'_i+\beta_ih)\\
=&\sum_{1\leq i\leq g}(\varphi(\overline e_i)+r_{2i-1})\beta_i\\
=&\sum_{1\leq i\leq g}r_{2i-1}\beta_i+\sum_{1\leq i\leq g}\beta_i\varphi(\overline e_i)\\
 =&r+\varphi(\overline t).
\end{array}$$

\item Let us prove i) simply for $t=0$ (since $K_0$ contains all $K_t$'s). By the propositions (\ref{prop:gs}) and (\ref{prop:genKt}), it suffices to prove that for any $Y\in\ker\pi_*$,  we have $(T_Y)_s(\ker\tilde\pi_*)\subset \ker\tilde\pi_*$.  This fact follows from the equality $(\ker\pi_*)^\perp=\ker\pi_*$ which induces that , for all $i$,
$(T_Y)_s(\overline e_i)=\overline e_i.$

\item In order to deduce ii) from i), we just have to notice that when $\varphi,\varphi'$ are in the same $G_s$-orbit, they necessarily have the same Arf-invariant, and to prove that when $\varphi,\varphi'\in{\cal E}_\pi$ have the same Arf-invariant, they are in the same $K_t$-orbit.

Assume that $\varphi,\varphi'\in{\cal E}_\pi$ and that $\alpha(\varphi'\circ s)=\alpha(\varphi\circ s)$, {\sl i.e.} by ii) $(\varphi'-\varphi)(\overline t)=0$. The difference between two quadratic forms is a linear fom, hence there exists a  vector $W$ defined by 
$$\forall x\in\Z_2^{2g}, (\varphi'-\varphi)(\overline x)=(\varphi'-\varphi)(s(x))=W.x$$
 hence $\varphi'\circ s=\varphi\circ s\circ T_W$ {\sl i.e.} $\varphi'=\varphi\circ (T_W)_s$ (see proof of the proposition 16 in [BGHM]). By definition of $W$, we have $W.t=(\varphi'-\varphi)(\overline t)=0$ and  for all $x\in \ker(\pi_*), W.x=(\varphi'-\varphi)(\overline x)=0$, hence $T_W\in K_t$.
\end{itemize}

\qed

\bs

\begin{cor}\label{cor:1,24}
Let $\psi$ and $\psi'$, be two elements of ${\cal E}(P_n).$ There exists $f\in G_s$  such that $\psi\circ\tilde\pi_*\circ f_s=\psi'\circ\tilde\pi_*$ if and only if $\sum_0^g\beta_j(\delta_0+\delta_j)=0$.

It is possible to choose $f$ to be a transvection $T_V\in K_t$ with $V=\sum_0^g(\delta_0+\delta_j)e'_j,$ (hence $\pi_*\circ f=\pi_*$).

\end{cor}

\subsection{Comparaison with the weak equivalence of coverings}\label{ssec:we}

\bs

\begin{nota}\label{nota:real}
Let $\tilde F$ be an automorphism of $H_1(P_n;\Z_2)$. We will say that $\tilde F$ is {\sl realizable} if there exists  $\Gamma$ an automorphism of $\pi_1P_n$ inducing $\tilde F$ in homology.
\end{nota}

\bs

\begin{thm}\label{lem:ader} Let  $\tilde F$ be an automorphism of $H_1(P_n;\Z_2)$, $F$ the endomorphism of $H_1(N_{g+1};\Z_2)$ and $\delta$ the linear form on $H_1(N_{g+1};\Z_2)$ defined by
$$\forall x\in H_1(N_{g+1};\Z_2),\quad \tilde F(\overline x)=\overline{F(x)}+\delta(x)h,$$
(for the notation $\overline x$, see Notation (\ref{nota:sigm})).
$\tilde F$ is {\rm realizable} if and only if  $\tilde F(h)=h$, $F$ is orthogonal, and $\delta(\sum_{i=0}^gv_i)=0$.
\end{thm}

\proof\quad 1) Let us prove that if  $F$ is orthogonal, $\tilde F(h)=h$ and $\delta(\sum_{i=0}^gv_i)=0$, then $\tilde F$ is realizable.
Because $F$ is orthogonal, there exists (\cite{GP}) an automorphism $\gamma$ of $\pi_1(N_{g+1})$ such that $\gamma_*=F$.
Let us choose integers $d_j\in\delta_j$ such that $\sum d_j=0\in\Z$ and  define $\Gamma_1(u_i)=u_ih^{d_i}$ and $\Gamma_1(h)=h$ (see the presentation (\ref{eqi:Pn})).  Now let $\Gamma_2$ be defined by $\Gamma_2=\gamma_*\times id_\Z$ and $\Gamma=\Gamma_2\circ\Gamma_1$. We obtain an automorphism $\Gamma$ of $\pi_1P_n$ which gives $\tilde F$ in homology.

2) Let us suppose now that  $\tilde F$ is  realizable and  $\Gamma$ is an automorphism of $\pi_1(P_n)$ such that $\Gamma_*=\tilde F$.

On order to prove the equality $\tilde F(h)=h$, we use the fact that $\Gamma$ lets invariant the center  $Z(\pi_1(P_n))$ of  $\pi_1(P_n)$. We have $Z(\pi_1(P_n))=Z(\pi_1(N_{g+1}))\times\Z$, and $Z(\pi_1(N_{g+1}))$ is easily computed after decomposing
$\pi_1(N_{g+1})$ as $F_g\star_\Z\Z$~:

if
 $g>1$, $Z(\pi_1(N_{g+1}))$ is trivial, hence $Z(\pi_1(P_n))$ is the sub-group generated by $h$, 

and if $g=1$, $Z(\pi_1(N_{g+1}))$ is the sub-group, isomorphic to $\Z$, generated by $x_0^2=x_1^{-2}$. Hence $Z(\pi_1(P_n))$ is the sub-group,  isomorphic to $\Z^2$, generated by $u_0^2=u_1^{-2}$ and $h$.

When $g>1$, it follows that $\Gamma(h)=h$ or $h^{-1}$, hence $\tilde F(h)=h$. 

When $g=1$, we only have  $\Gamma(h)=u_0^{2u}h^v$ (with $u,v\in\Z$), but moreover $h$ is not a square in $\pi_1(P_n)$, nor is $\Gamma(h)$, so $v$ is odd, and again $\tilde F(h)=h$.

Now we want to show that $\delta(\sum_{i=0}^gv_i)=0$. For all $\psi\in{\cal E}(P_n),$ we have (because $\tilde F(h)=h$) that $\psi\circ\tilde F\in{\cal E}(P_n)$, then by the easy part of theorem (\ref{thm:an}) below, $\psi\circ(\tilde F-id)(\sum_{i=0}^g\overline{v_i})=0$. From this follows that $0=(\tilde F-id)(\sum_{i=0}^g\overline{v_i})=\overline{(F-id)(\sum_{i=0}^g{v_i})}+\delta(\sum_{i=0}^g{v_i})h=0$, in particular $\delta(\sum_{i=0}^g{v_i})=0$. (An other proof is just to remark that if we denote  $\Gamma(x_i)=y_ih^{d_i}$ with $y_i\in\pi_1(N_{g+1})$, then we get $2\sum d_i=0$).

It remains to prove that $F$ is orthogonal. We have $F=\gamma_*$, where $\gamma$ is the endomorphism of $\pi_1(N_{g+1})$ deduced from $\Gamma$ after composition on each side by the natural projection and inclusion using that $\pi_1(P_n)=\pi_1(N_{g+1})\times\Z$. When $g>1$,  $\Gamma(h)$  has no component on $\pi_1(N_{g+1})$ (see above). Let us denote by $\gamma'$ the endomorphism of $\pi_1(N_{g+1})$ deduced from $\Gamma^{-1}$, we find the equality $\gamma\circ\gamma'=\gamma'\circ\gamma=id$, so $\gamma$ is an  automorphism, hence  $F=\gamma_*$ is orthogonal. In the particular case when $g=1$, since $F$ is injective and $(F-id)(v_0+v_1)=0$ (see above), $F$ is orthogonal.

\bs

We recall the classical definition of {\sl weak equivalence} between two coverings.
\begin{defi}\label{defi:wec}
Two coverings $E_1\stackrel{p_1}\longrightarrow P$ and $E_2\stackrel{p_2}\longrightarrow P$ are {\rm weakly equivalent} if there exist  an automorphism $\Gamma$ of 
$\pi_1 P$ and an isomorphism $\gamma$ between $\pi_1E_1$ and $\pi_1E_2 $ such that  $p_2\circ\gamma=\Gamma\circ p_1.$ 
\end{defi}

 In our particular case where the  two   $\Z_2-$coverings are special, they
correspond to morphisms $\psi$ and $\psi'\colon \pi_1P_n\to \Z_2$. They
are
 wealyk equivalent if and only if there exists  an automorphism $\Gamma$ of
 $\pi_1 P_n$ such that $\psi\circ \Gamma=\psi'$. Using the same letter
$\psi$
 and $\psi'$ for the corresponding morphisms $H_1(P_n;\Z_2)\to \Z_2$, the
two
 special $\Z_2-$coverings are  weakly equivalent if and only there exists a
 realizable automorphism $G$ of $H_1(P_n;\Z_2)$ such that $\psi\circ
 G=\psi'$.

\bs

\begin{thm}\label{thm:an}
The following four properties are equivalent  :
\begin{description}
\item[1)]  Two  special  $\Z_2-$coverings $E_{\psi}$ and $E_{\psi'}$ are weakly equivalent;
\item[2)]  Two  special  $\Z_2-$coverings $E_{\psi}$ and $E_{\psi'}$ are $\pi_1-$isomorphic;
\item[3)]  $\Sigma_0^g \psi(\overline v_i)=\Sigma_0^g \psi'(\overline v_i)$;
\item[4)]  There exist choices of the linear lifting $\sigma(\psi,\psi')$ for which $E_{\psi}$ and $E_{\psi'}$ are in the same class under the action $A_1^{\sigma(\psi,\psi')},$ (see Definition (\ref{defi:1}).
\end{description}
\end{thm}

\proof 

$1)\Rightarrow 2)$ is evident and $2)\Leftrightarrow 3)$ follows from the theorem (\ref{thm:kernon}).

We will show that $3)\Rightarrow 1).$ Let $\psi$ and $\psi'$ be two elements of ${\cal E}(P_n)$ and  $\delta_j=(\psi'-\psi)(\overline v_j)$. We consider   the automorphism $\tilde F$ of $H_1(P_n;\Z_2)$ with matrix $\left( \begin{array}{cc}
                                  Id&0\\
                               \delta & 1
                                 \end{array}
                                  \right)$ in the basis $(\overline v_0,\ldots,\overline v_g,h)$ then the theorem (\ref{lem:ader}) shows that $\tilde F$ is realizable.

$3)\Leftrightarrow 4)$ follows from the theorem (\ref{thm:1}). More precisely, $4)\Rightarrow 3)$ is straightforward, and to prove the converse, we just have to show that if $3)$ holds and if $\psi'\neq\psi$, there exists a $\sigma$ for which the associated $\psi_0,\psi_1$ differ from $\psi$ and $\psi'$. By assumption, there exist at least two indices $i$ such that $\psi'(\overline v_i)\neq\psi(\overline v_i)$, for instance $i=0$ and $i=1$. Choose $\rho_0=\psi(\overline v_0)$, $\rho_1=\psi'(\overline v_1)$, $\rho_2,\ldots,\rho_g$ arbitrary, and take the corresponding $\sigma$ (see Notation (\ref{nota:sigm})).

\qed

\bs

\subsection{Relation between $ Sp(\Z_2,2g)$ and $O(\Z_2^{g+1})$}\label{ssec:birm}

\bs

\begin{defi}\label{defi:lift}
When  an endomorphism $f\colon H_1(F_g;\Z_2)\to  H_1(F_g;\Z_2)$ and an endomorphism $F\colon H_1(N_{g+1};\Z_2)\to  H_1(N_{g+1};\Z_2)$  are related by  $F\circ \pi_\star=\pi_\star\circ f$, we will say that {\rm $f$ is a lifting of $F$} or that {\rm $F$ is a projection of $f$}
\end{defi}

\bs

\begin{thm}\label{thm:symsym} 

The basis $(e_i,e'_i)$ of $H_1(F_g;\Z_2)$ and $(v_0,v_0+v_i)$ of $ H_1(N_{g+1};\Z_2)$ are chosen as above.
\begin{description}
\item[1)] An endomorphism $f\colon H_1(F_g;\Z_2)\to  H_1(F_g;\Z_2)$ and an endomorphism $F\colon H_1(N_{g+1};\Z_2)\to  H_1(N_{g+1};\Z_2)$  are related as in the above definition if and only if, their matrices are respectively $M_f=\left( \begin{array}{cc}
 A &0 \\ 
 C &D
\end{array}
 \right)$ and $M_F=\left( \begin{array}{cc}
 b &0 \\
 B &A
\end{array}
 \right)$ with the same   $g\times g$ matrix $A$ with coefficients in
$\Z_2$.
\item[2)]  Let $f$ be an endomorphism of $H_1(F_g;\Z_2)$ with matrix,  $M_f=\left( \begin{array}{cc}
 A &0 \\
 C &D
\end{array}
 \right), $ and $A$ invertible. $f$ admits an orthogonal projection $F$ if and only if $A^tSA=S, $ where the matrix $S$ is obtained by restricting to $Im(\pi_\star)$ the scalar product defined on $H_1(N_{g+1};\Z_2)$.
\item[2)bis] There exists $F\in O(\Z_2^{g+1})$, extending  the automorphism ${\hat F}$ of $Im (\pi_\star)=H$, found in the part 2) of the theorem (\ref{thm:symsym}), unique if $g$ is even and double if $g$ is odd
\item[3)]    Let $F$ be an endomorphism of $ H_1(N_{g+1};\Z_2)$. $F$ admits a symplectic lifting $f$ if and only if $F(Im(\pi_\star))=Im(\pi_\star).$  
\end{description}
\end{thm}

\bs

\begin{rem}\label{rem:cite}
\begin{itemize}
\item If $f$ is a symplectic automorphism of $H_1(F_g;\Z_2)$,  then its matrix $M_f=\left( \begin{array}{cc}
 A &0 \\
 C &D
\end{array}
 \right), $ verifies $M^t_f.K.M_f=K$ where $K$ is the matrix of the symplectic product of $H_1(F_g;\Z_2)$. The  translation on   $A,0,C,D$ of this condition is :
$D^t.A=I$ and $C^t.A=A^t.C$. There exists such $C,D$ if and only if $A$ is invertible.
\item If $F$ is an orthogonal automorphism of $ H_1(N_{g+1};\Z_2)$ then $F$ is bijective and fixes the element $v_0+\cdots+v_g$ (see \ref{prop:mi}) i.e. lets stable  $(v_0+\cdots+v_g)^\bot=Im(\pi_\star)$. 
\item Let $C$ be a $g\times g$ symetric matrix, all the symplectic automorphims of $H_1(F_g;\Z_2)$ with matrix $M_f=\left( \begin{array}{cc}
 Id &0 \\
 C & Id
\end{array}
 \right), $   are liftings of the identity of $ H_1(N_{g+1};\Z_2)$.
\end{itemize}
\end{rem}

\bs

{\sl Proof
of  theorem} (\ref{thm:symsym}).

\begin{description}
\item[1)] $f(\ker\pi_\star
)\subset\ker\pi_\star$ and $F(Im\pi_\star )\subset Im\pi_\star $ if and only if the matrices of $f$ and 
$F$ are respectively of the form $M_f=\left( \begin{array}{cc}
 A &0 \\ 
 C &D
\end{array}
 \right)$ and $M_F=\left( \begin{array}{cc}
 b &0 \\
 B &A'
\end{array}
 \right).$
We have $\pi_\star
f(e'_j)=0=F\pi_\star (e'_j)$, $\pi_\star
f(e_j)=\pi_\star (\sum_i A_{i,j}e_i+C_{i,j}e'_i)=\sum_iA_{i,j}(v_0+v_i)$ and
$F\pi_\star (e_j)=F(v_0+v_j)=\sum_i A'_{i,j}(v_0+v_i)$ hence $\pi_\star f=F\pi_\star $ if and only if
$A'=A$. 
\noindent\item[2)]
 i) 
Let us suppose that there exists an orthogonal automorphism of $ H_1(N_{g+1};\Z_2)$ such that $F\pi_\star =\pi_\star f.$
In the basis $v_0+v_i=\pi_\star(e_i)$ of $Im(\pi_\star)$,  the matrix $S$, obtained by restricting the scalar product to $Im(\pi_\star)$, has entries $S_{ij}=(v_0+v_i)(v_0+v_j)=1+\delta_{ij}.$ 
The automorphism $F$, restricted to $Im(\pi_\star)$ respects the scalar product and  $F_{\mid Im(\pi_\star)}=A$ so $A^tSA=S.$

\noindent\noindent ii) Let us suppose that  $A\in  Gl(\Z_2^g)$ and $A^tSA=S. $  This matrix $A$ can be considered as the matrix of an automorphism $\hat F$ of $Im(\pi_\star)$, respecting the scalar product. The following general lemma (\ref{lem:gene}) will prove that $\hat F$ can be extended as an orthogonal automorphism $F$ such that $F\circ \pi_\star=\pi_\star\circ f$.
\noindent\item[3)]
 With a $g\times g$ matrix $A$  invertible, now we know how to construct a symplectic automorphism $f$ such that $F\pi_\star =\pi_\star f$. The matrices $\left( \begin{array}{cc}
 b &0 \\
 B &A
\end{array}
 \right),$ with $A$ invertible are exactly the matrices of the $F$ such that $F(Im\pi_*)= Im\pi_*.$
\end{description}
\qed

\bs

 The orthogonal  of $Im(\pi_\star)$ is $Im(\pi_\star)^\bot=\Z_2E,$ with $E=\Sigma_{i=0}^gv_i=(g+1)v_0+\pi_\star (e)$ for $ e=\Sigma_{i=1}^ge_i.$  Let us remark here that the generator $E$ of $Im(\pi_\star)^\bot $ belongs to $Im(\pi_\star)$ if and only if the genus $g$ is odd. 

\bs

\begin{lem}\label{lem:gene}  In $\Z_2^{g+1}$ equiped with a scalar product $x.y=\Sigma x_iy_i$, let $H=E^\bot$ be any hyperplane and $\hat F\in O(H)$. For any $v\notin H,$ the  following is true:
\begin{enumerate}
\item[ i)] $\exists v'\in \Z_2^{g+1},\forall x\in H, v'.{\hat F}(x)=v.x;$
\item[ ii)] there are two solutions $v'$ differing by $E$. If $E\notin H,$  only one of them verifies $v'\notin H.$ If $E\in H$ then  both verify $v'\notin H.$
\item[ iii)] Moreover, in the particular case where $E=(1,\cdots ,1)$ one (or the two) preceeding solutions verify $v'.v'=v.v$.
\end{enumerate}
\end{lem}

\proof 
\begin{enumerate}
\item[ i)] If $E\notin H, v':={\hat F}(v-E)$ is a solution. If $E\in H$ let $K=H\cap v^\bot$, then $H=K\oplus \Z_2E.$ We have ${\hat F}(E)\notin {\hat F}(K)$ then ${\hat F}(K)^\bot \not\subset {\hat F}(E)^\bot$ hence there exists $v'\in {\hat F}(K)^\bot, v'\notin  {\hat F}(E)^\bot$ i.e. $\exists v',\forall x \in K\cup \{E\}, v'.{\hat F}(x)=v.x.$
\item[ ii)] If $v'_0$ is a solution, then $v'$ is also a solution if and only if $v'-v'_0\in Im({\hat F})\bot=H^\bot=\Z_2E$, it means $v'=v'_0$ or $v'=v'_0+E.$ When $E\notin H, (v'+E).E=(v'.E)+1$ then only one of them verifies $v'.E=1.$ When $E\in H$, we have ${\hat F}(E)=E$ then both of them verify $v'.E=v'.{\hat F}(E)=v.E=1.$
\item[ iii)] When $E=(1,\cdots ,1), \forall x\in \Z_2^{g+1},x.x=x.E$ this implies that $\forall v,v'\notin H, v'.v'=v'.E=1=v.E=v.v.$
\end{enumerate}
\qed

This lemma proves that there exists $F\in O(\Z_2^{g+1})$, extending  the automorphism ${\hat F}$ of $Im (\pi_\star)=H$, found in the part 2) of the theorem (\ref{thm:symsym}), unique if $g$ is even and double if $g$ is odd and completes the proof of  theorem (\ref{thm:symsym}).

\bs

What follows may be applied to the endomorphism $\tilde f=f_s$ of $H_1(P_o;\Z_2)$, with symplectic $f$, but whithout any simplification.

\bs

\begin{prop}\label{prop:reali} 
\begin{enumerate}
\item[ a)] An endomorphism $\tilde f$ of $H_1(P_o;\Z_2)$ verifies $\tilde f(ker(\tilde\pi_*))\subset ker(\tilde\pi_*)$ if and only if there exists an endomorphism $\tilde F$ of $H_1(P_n;\Z_2)$ such that $\tilde\pi_*\circ\tilde f=\tilde F\circ\tilde\pi_*$.
\item[ b)] Such $\tilde f$ which moreover fixes $h$ has its  matrix in the basis $(\overline{e_i}_{1\leq i\leq g},\overline{e'_i}_{1\leq i\leq g},h)$ of the following form
$\left( \begin{array}{ccc}
 A &0&0 \\ 
 C &D&0\\
k&0&1
\end{array}
 \right).$
 In the basis $(\overline{v_0},\overline{v_0+v_i}_{1\leq i\leq g},h),$ the matrices of the associated $\tilde F,$ are  $\left( \begin{array}{ccc}
 b &0&0 \\
 B &A&0\\
\lambda&k&1
\end{array}
 \right),$ 
with arbitrary $b,B,\lambda.$ 
\end{enumerate}\end{prop}

\proof\quad $a)$ is evident and the proof of $b)$ is the same as the proof of the item $1)$ of  theorem (\ref{thm:symsym}).
\qed

\bs

\begin{thm}\label{thm:der} Let $\tilde f$ be like in the $a)$ and $b)$ of the above proposition. The conditons to get that one of its associated $\tilde F$ is realizable are : the invertibility of the matrix $A$ with $A^t.S.A=S,$ where $S$ is the restriction to $Im(\pi_\star)$ of the intersection matrix on $H_1(N_{g+1};\Z_2)$, like in the theorem (\ref{thm:symsym}), and if $g$ is odd, $\sum_{i=1}^g k_i=0$. The number of realizable $\tilde F$ is $1$ if $g$ is even and  $4$ if $g$ is odd.
\end{thm}

\proof\quad

Given $\tilde f$ (in particular $A$ and $k$), let us denote by $F$ the endomorphism (depending of the choice of $b,B$) of $H_1(N_{g+1};\Z_2)$ with matrix  $\left( \begin{array}{cc}
 b &0 \\
 B &A
\end{array}
 \right),$ in the basis $(v_0,(v_0+v_i)_{1\leq i\leq g})$ and $\delta$ (depending of the choice of $\lambda$) the linear form on $H_1(N_{g+1};\Z_2)$ defined by $\delta(v_0)=\lambda$, $\delta(v_0+v_i)=k_i$. From  theorem (\ref{lem:ader})
$\tilde F$ is realizable if and only if $F$ is orthogonal (which is a constraint on $A$ and on the choice of $b,B$) and 
$0=\delta(\sum_{i=0}^g v_i)=(1-g)\lambda+\sum_{i=1}^g k_i$ (a constraint in the choice of $\lambda$). 

From the item $2)bis$ of theorem (\ref{thm:symsym}), the number of  solutions $(b,B)$ is
\begin{itemize}
\item  $0$ if the  condition ``$A$ invertible and  $A^t.S.A=S$'' is not true, 

\item $1$  if this condition is true and $g$ is even,

\item $2$  if this condition is true and  $g$ odd. 
\end{itemize}

The number of solutions $\lambda\in\{0,1\}$ is $1$ if $g$ is even, $0$ if $g$ is odd and $\sum_{i=1}^g k_i=1$, $2$ if $g$ is odd and $\sum_{i=1}^g k_i=0$.

\qed

\bs

{\bf Application of the above study to the weak-equivalence of special  $\Z_2-$coverings.}

\bs

\begin{thm}\label{thm:fin}
\begin{itemize}
\item If $g$ is even, for each $\psi\in {\cal E}(P_n)$, the two $\psi'$ such that $\psi\circ\tilde\pi_*\circ\tilde f=\psi'\circ\tilde\pi$ (see proposition (\ref{prop:2g})) are not in the same  weak-equivalence classes, hence $\psi$ is weak-equivalent to exactly one of them, and (provided that $A$ is invertible and  $A^t.S.A=S$), this weak-equivalence can be explicited from $\tilde f$.
\item If $g$ is odd, the weak-equivalence  class of $\psi$  depends only on $\psi\circ\tilde\pi$ hence for each 
 $\psi$, the two $\psi'$ such that $\psi\circ\tilde\pi_*\circ\tilde f=\psi'\circ\tilde\pi_*$ are either together in the same weak-equivalence  class of $\psi$, or together in the other one. In the first case, for each of them, it is possible (provided that $A$ is invertible and $A^t.S.A=S$) to build from $\tilde f$ two equivalences between $\psi$ and $\psi'$. 
\end{itemize}
\end{thm}



\begin{thebibliography}{VBGT}

\bibitem[Arf]{Arf} Arf C. :Untersuchungen \:uber quadratische Formen in Korpern der charakteristik 2.I  J. Reine Angew. Math. 183 (1941) 148-167
\bibitem[A]{A} Atiyah M.: Riemann surfaces and Spin
      structures. Ann. scient. Ec. Norm. Sup., 4 s\'erie, t. 4 (1971) 47-62
\bibitem[BGHM]{BGHM} Bauval A. ; Gon\c calves D. L. ; Hayat C. ; Mello M. H. P. L. :A classification of special 2-fold coverings. Geometry and Topology Monographs 14 (2008) 27-47
    URL: http://www.msp.warwick.ac.uk/gtm/2008/14/p002.xhtml
\bibitem[BC]{BC} Birman J. S. ; Chillingworth D. R. J. : On the homeotopy group of a non-orientable surface. Proc. Camb. Phil. Soc. 71 (1972) 437-448
\bibitem[GP]{GP} Gadgil S. ; Pancholi D. : Homeomorphisms and the homology of non-orientable surfaces. 
Proc. Indian Acad. Sci. (Math. Sci.) Vol. 115, No. 3 (2005) 251-257
\bibitem[GHM]{DCM} Gon\c calves, D.L.: Hayat C.; Mello, M.H.P.L.: Spin-structures and 2-fold coverings, BSPM, 23 1, 2 (2005) 29-40
\bibitem[J1]{Joh1} Johnson, D.: Homeomorphisms of a surface which act trivially on homology. Proc. of the AMS. 75 1 (1979)119-125
\bibitem[J2]{Joh2} Johnson, D.: Spin-structures and quadratic forms on
      surfaces. J. London Math. Soc. (2) 22 (1980) 365-373
\bibitem [K1]{K1} Korkmaz M. : First homology group of mapping class groups of nonorientable surfaces. Math. Proc. Camb. Phil. Soc.123 (1988) 487-499
\bibitem[K2]{K2} Korkmaz M. : Mapping class groups of nonorientable surfaces. Geometrica Dedicata 89 (2002) 109-133
\bibitem[LI]{LI} Lickorish W. B. R. On the homeomorphisms of a non-orientable surface. Proc. Camb. Phil. Soc. 61 (1965) 61-64
\bibitem[MCP]{MCP} McCarthy J. D. ; Pinkall U. : Representing homology automorphisms of nonorientable surfaces. Max Planck Institut preprint MPI/SFB 85-11
\bibitem[Mi]{Mil} Milnor, J.: Spin structure of manifolds.  Enseignement des math.
2 9 (1963) 198-203 
\bibitem[ZVC]{ZVC} Zieschang, H.; Vogt, E.; Coldewey, H.-D.: Surfaces and planar
      discontinuous groups. Lecture Notes Math. 835,
      Berlin-Heidelberg-New York: Springer-Verlag 1980
\end{thebibliography}
\end{document}